\documentclass[12pt]{article}
\usepackage{subfigure}
\usepackage{algorithmic}
\usepackage{graphicx} 
\usepackage{amsmath}
\usepackage{amssymb}
\usepackage{amsfonts}
\usepackage{amsthm}
\usepackage{algorithm}
\usepackage{algorithmic}
\usepackage{authblk}
\usepackage{hyperref}
\newtheorem{theorem}{Theorem}[section]

\newcommand{\no}{\nonumber}
\newcommand{\bb}{\boldsymbol}
\newtheorem{exmp}{Example}

\begin{document}
\title{Statistical test for detecting community structure in real-valued edge-weighted graphs}
\author{Tomoki Tokuda \thanks{tomoki.tokuda@oist.jp}}
\affil{Okinawa Institute of Science and Technology\\ Graduate University, JAPAN}

\date{October 11, 2016}
\maketitle

\begin{abstract}
We propose a novel method to test the existence of community structure 
of undirected real-valued edge-weighted graph. The method is based on Wigner semicircular law on
the asymptotic behavior of the random distribution for eigenvalues of a real symmetric matrix. 
We provide a theoretical foundation for this method and report on its performance 
in synthetic and real data, suggesting that our method outperforms other state-of-the-art methods. 
\end{abstract}

\section{Introduction}
Clustering objects based on their similarities is a basic data mining approach in statistical analysis. In particular, graph data (or, network data) that reflect relationships between nodes, are often acquired in various scientific domains such as protein-protein interaction, neural network and social network \cite{fortunato2010community}, which potentially provides quite useful information on the underlying structure of the system in question.

\begin{figure}[ht!]
     \begin{center}
        \subfigure[Edge-weighted graph]{
            \label{graphrep}
               \includegraphics[scale=0.16, trim=0mm 0mm 0mm 0mm]{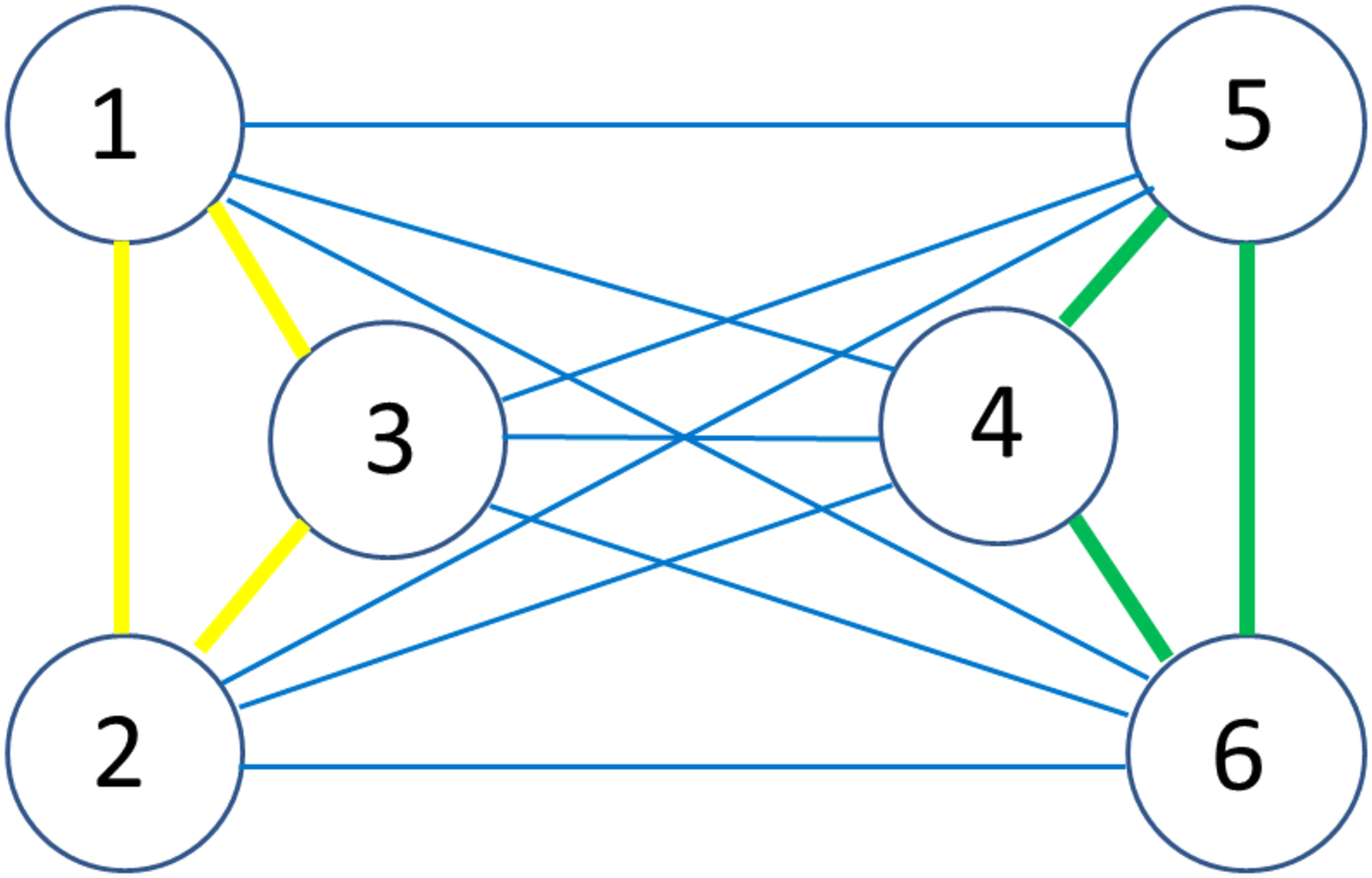} 
        }
        \subfigure[Edge-weight matrix]{
           \label{matrixrep}
           \includegraphics[scale=0.12, trim=0mm 0mm 0mm 0mm]{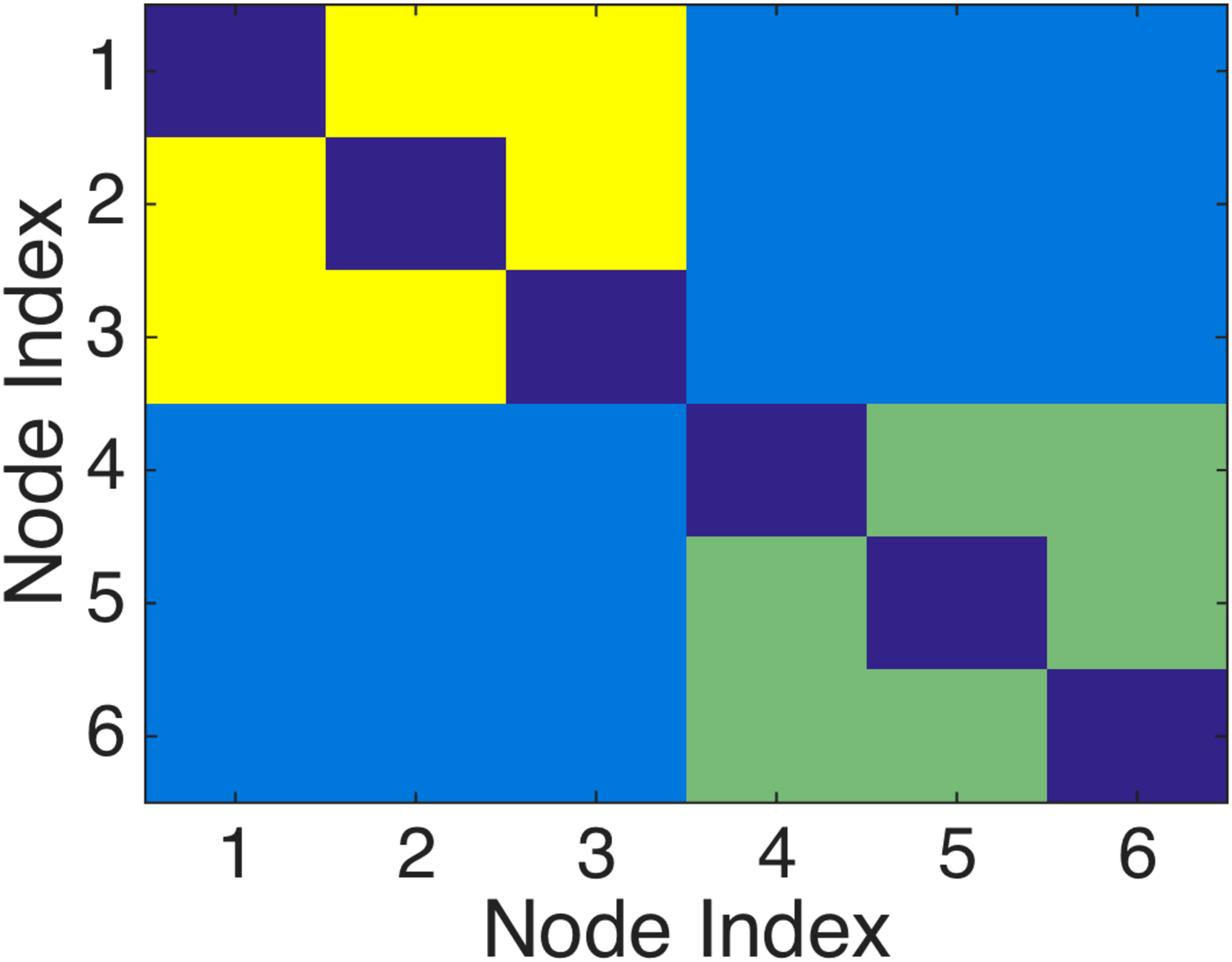} 
        }\\
    \end{center}
    \vspace{-5mm}
        \caption{\it \small Illustration of two-way community structure in a graph. Panel (a): Graph representation 
        (edge-weighted graph). Panel (b): Matrix representation (edge-weight matrix), where strengths of relationships between nodes are denoted by color.}
   \label{repmat}
\end{figure}

Specifically, our interest is to detect a possible \lq community\rq,~or cluster structure of undirected graph, which is defined as block structure of a graph (Fig.\ref{repmat}a), where the corresponding edge-weight matrix consists of several cluster blocks (four cluster blocks in Fig.\ref{repmat}b).  To detect such structure, a number of clustering methods have been proposed in the literature of statistical physics and information theory
\cite{newman2004finding,  reichardt2006statistical, bolla2013spectral}. Mainly, there are four approaches: graph partitioning, hierarchical clustering, partitional clustering and spectral clustering \cite{fortunato2010community, bolla2013spectral, ng2002spectral}. 

However, the conventional framework for analysis of community structure is typically for an unsigned graph in which 
an edge weight is constrained to be non-negative. Recently, it has gained much attention to
analyze a signed graph that allows for negative weights  \cite{kunegis2010spectral}. Indeed, in real data, it is often essential  to take into account  negative as well as positive relationship for a better understanding of the underlying community structure in a graph such as social network. Most methods in literature, however, address this problem in a rather limited framework in which edge weights within a cluster are positive while those between clusters negative (i.e., weakly balanced structure) \cite{kunegis2010spectral}. On the other hand, it still remains an open question as how to cluster nodes in a more general framework such as negative edge weights within a cluster \cite{yang2007community}. 

In the present paper, we consider a general framework for community structure as follows. We assume that edge weights are independently generated from a generative model that is specific to a particular cluster block, which characterizes a distribution of edges in each cluster block. Further, we assume that these distributions are distinguishable in terms of their mean and variance. 
For this framework, as a first step of addressing a clustering problem, we aim to develop a statistical method for testing the existence of the underlying community structure. 

As regards statistical test on community structure, several methods have been proposed in the context of unsigned (weighted or unweighted) graph \cite{fortunato2010community}. A major approach to this problem is to evaluate the stability of cluster solutions when the graph in question is contaminated with noise \cite{gfeller2005finding, karrer2007robustness}. 
If similar cluster solutions are yielded for contaminated graphs, it suggests the stability of the cluster solution for the original graph, providing the evidence of the community structure. 
The bootstrap method \cite{rosvall2010mapping} is in the similar line with this approach. A second approach is based on comparison of the cluster solution for the original graph with those solutions of randomly permuted graphs. 
As a statistic for testing the significance, the entropy of graph configurations \cite{bianconi2009assessing}, or
\lq C-score\rq~focusing on the lowest internal degrees \cite{lancichinetti2010statistical} have been proposed.
The common feature of these state-of-the-art methods is that a cluster solution to the graph in question is required for testing. In other words, the test result depends on a clustering method that one uses. In this sense, these methods test the significance of a yielded cluster solution, rather than the existence of community structure itself. For the general framework of our interest, such an approach is not applicable because appropriate clustering methods are not readily available.  

We propose a general method for testing community structure of edge-weighted graph with real-valued weights, which does not require a cluster solution. Our method is based on the asymptotic behavior of eigenvalues of the normalized weight matrix of graph, which is described by Wigner semicircular law when there is no community structure. 
In the similar line with our approach, in the case of binary-valued graph, a statistical test on community structure has been recently proposed \cite{bickel2016hypothesis} that is based on the exact asymptotic behavior of (maximum) eigenvalues.
However, their method is not directly applicable to real-valued graphs where we take into account
both mean and variance, because Bernoulli distribution assumed in their method cannot  properly capture these quantities.
Our method provides a nontrivial extension of detecting community structure to real-valued graphs, which has a wide range of applications to network data. 
In the following sections, first, a theoretical foundation for our method is provided. Second, it is shown that our method outperforms other methods in synthetic data. Third, we apply our method to real data. 

\section{Method}
Our statistical test on community structure is based on the probability distribution of eigenvalues of the normalized edge-weight matrix
(we define \lq normalization\rq~in Section~\ref{stattest}). 
We make the best use of asymptotic results on such a distribution when there is no community structure, which have been intensively studied in the field of Random Matrix Theory of Theoretical Physics \cite{mehta2004random}. In this section, we provide a theoretical foundation for our statistical test.
\begin{figure}[ht!]
     \begin{center}
        \subfigure[Setting of parameters]{
            \label{lattice}
           \includegraphics[scale=0.16, trim=-0mm 0mm 0mm 0mm]{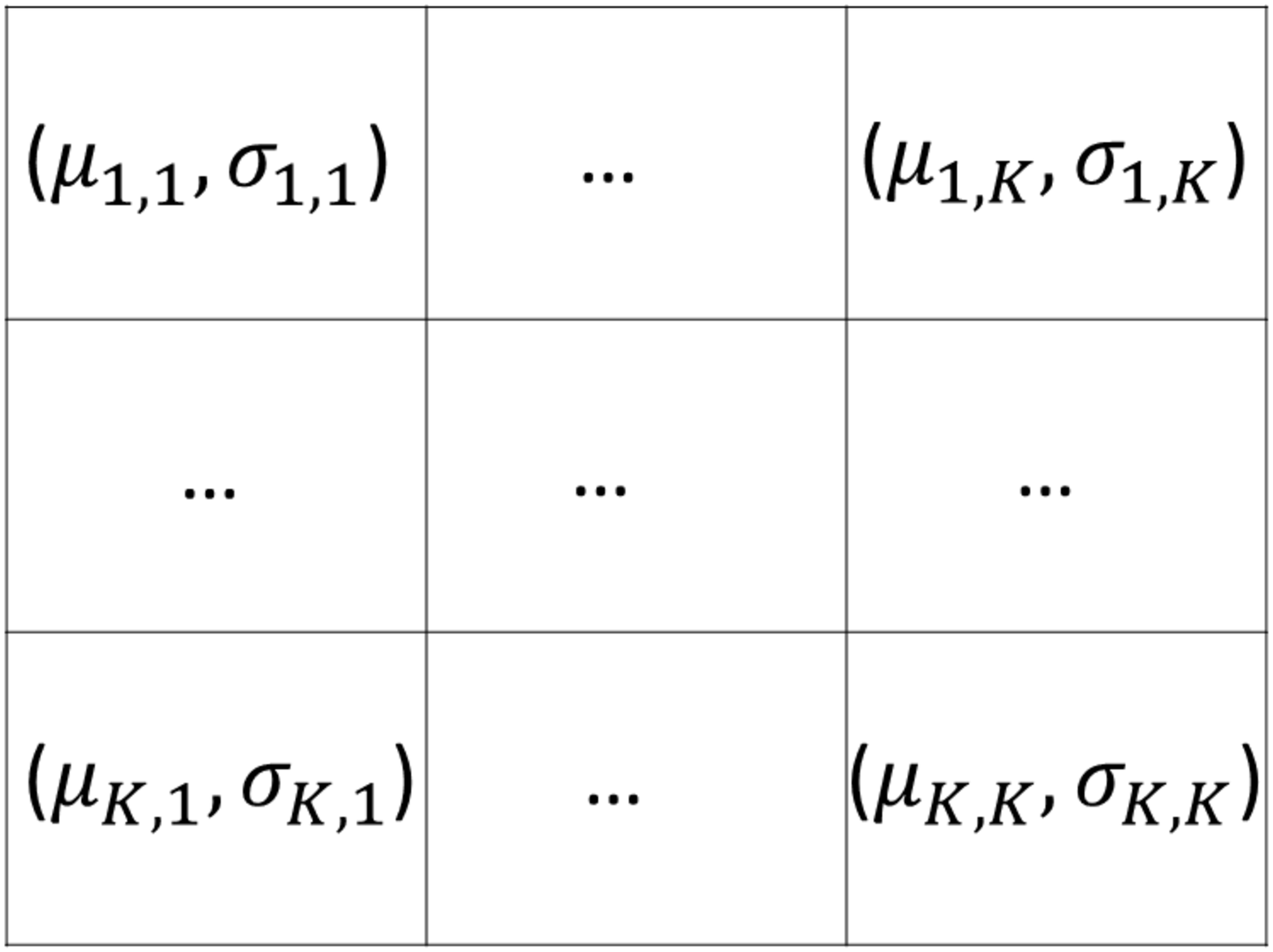} 
        }
        \subfigure[Tracy-Widom distibution]{
           \label{density}
           \includegraphics[scale=0.12, trim=0mm 0mm 0mm 0mm]{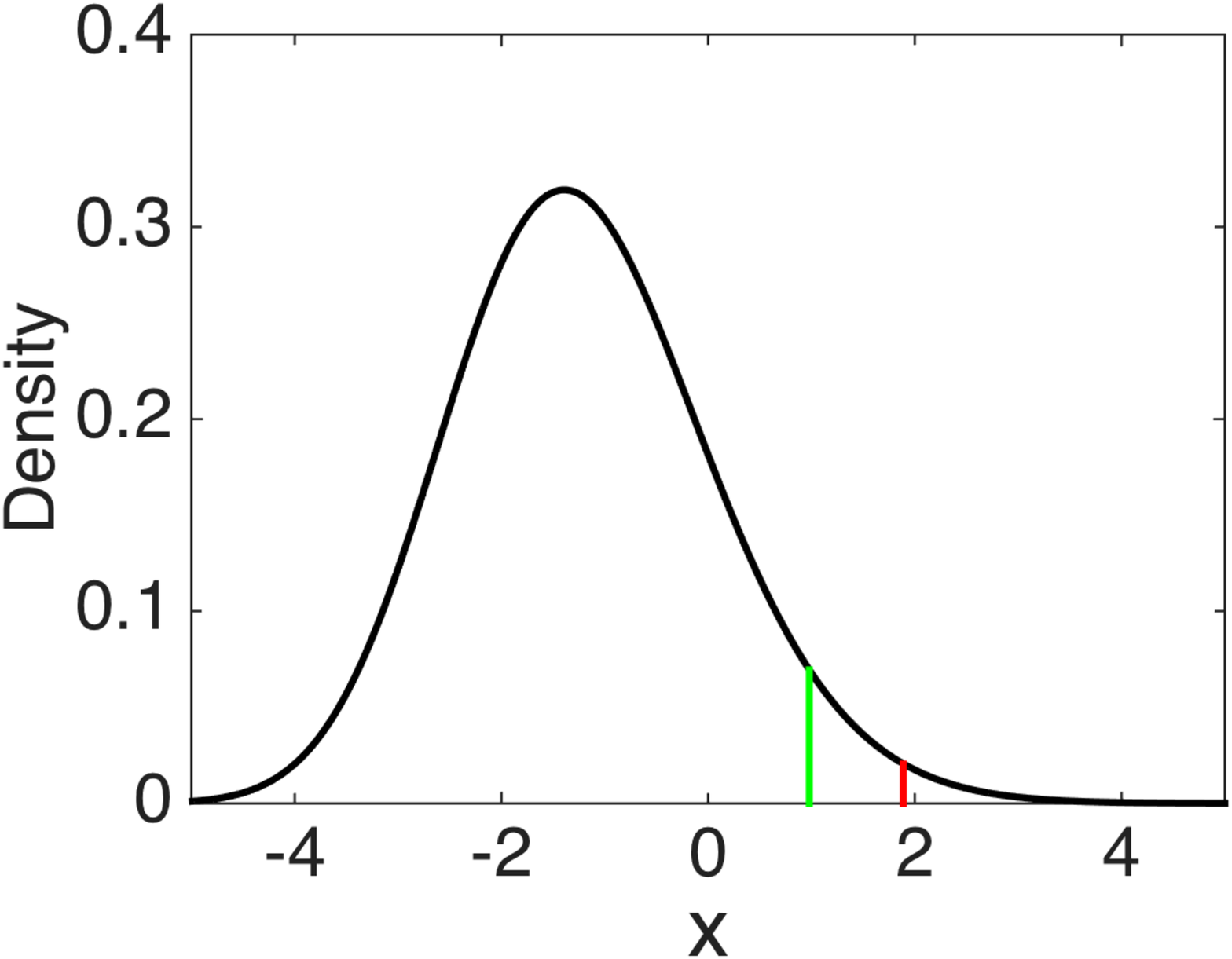} 
        }\\ 
    \end{center}
    \vspace{-5mm}
        \caption{\it \small Panel (a) : Illustration of setting of community structure in matrix representation
        where the nodes are arranged in the order of cluster labels. Each cluster block is characterized by 
        mean $\mu$ and standard deviation $\sigma$ with cluster block index $(k, k')$. 
        Panel (b): The density function of Tracy-Widom distribution for Gaussian orthogonal ensembles with $\beta=1$ (the first derivative of $F_1(x)$ in Eq.(\ref{F1})), generated by the function {\it dtw} in R-package \{RMTstat\}. The critical values at significance level $\alpha=0.05$ and 
        $\alpha/4 = 0.0125$ are $0.979$ (denoted in green line) and $1.889$ (in red line), respectively. 
          }
\label{F2}
\end{figure}

\subsection{Setting} \label{setting}
We consider a clustering problem of nodes for undirected edge-weighted graph $G=(V, E)$ where $V$ consists of $n$ vertices $\{v_1, \ldots, v_n\}$, and $E$ is represented by the edge-weight matrix $\bb{W}_n$, which is a $n \times n $ symmetric (real Hermitian) matrix with elements $w_{i, j}=w_{j, i} \in \mathbb{R}$ and $w_{i, i} = 0$ ($\mathbb{R}$ denotes a set of real numbers). Let us assume that there are 
$K$ clusters of nodes, denoting them as $c_1, \ldots, c_K$. We define a cluster block $(k, k')$ as a set of weights $w_{i, j}$ such that nodes $i$ and $j$ belong to the cluster $c_k$ and $c_{k'}$, respectively: $v_i \in c_k$ and $v_j \in c_{k'}$ ($1 \leq k, k' \leq K$). 
Here, we assume that each off-diagonal weight $w_{i, j}$ is independently drawn from a certain distribution. With this assumption, we define a $K$-way community structure as characterized by different distributions in $K \times K$ cluster blocks.
To elaborate this definition, we assume the following distribution for each cluster block:
\begin{eqnarray}
 \no w_{i, j} &\sim& g_{k, k'} ~~ (i \neq j)\\
  g_{k, k'} &=& \mu_{k, k'} + g \times \sigma_{k, k'},
  \label{settingw}
\end{eqnarray}
where  $v_i \in c_k$, $v_j \in c_{k'}$, and $g$ is a certain probability distribution.
This definition suggests that a pair of parameters $(\mu_{k, k'}, \sigma_{k, k'})$ characterizes each cluster block, hence, 
community structure (Fig.\ref{F2}a). 
Note that in this definition we exclude the degenerate case that 
$\mu_{k, k'} =constant$ and $\sigma_{k, k'}=0$ in which the variances become zero for the whole set of $\{w_{i, j}\}$. 

Since the community structure of our interest is based on differences of weight distributions, it is translation and scale  invariant for the whole weights. Hence, to simplify the problem, as a preprocess, we standardize off-diagonal elements of $\bb{W}_n$ using all off-diagonal weights $w_{i, j}(i \neq j)$ so that the mean is zero and the variance one. We denote 
as $S$ the mapping that standardizes edge-weight matrix in this way, transforming each element of the matrix as
\begin{eqnarray}
 \no  S:  &&w_{i, j} \rightarrow (w_{i, j} - \mu)/\sigma ~~\mbox{for}~ i \neq j \\
            &&w_{i, i}  \rightarrow 0,
 \label{S}
\end{eqnarray} 
where $\mu$ and $\sigma$ are the mean and the standard deviation of the whole off-diagonal elements $\{ w_{i, j }\}$. Practically, these mean and standard deviation may be replaced by the empirical counterparts
$\mu_{emp}$ and $\sigma_{emp}$. 
For the standardized edge-weight matrix $S(\bb{W}_n)$, we assume that the mean and the variance of $g$ in Eq.(\ref{settingw}) are zero and one, respectively. In this setting, the mean and the variance in cluster block $(k, k')$ are $\mu_{k, k'}$ and $\sigma_{k, k'}^2$, respectively. The differences of these parameters 
distinguish between clusters 
in terms of the first and second moments, while controlling higher 
moments than two. 
Using this setting of community structure, we define no community case as a single community with $K=1$ where  $\mu_{k, k'}=0$ and $\sigma_{k, k'}=1$ for $S(\bb{W}_n)$. Note that since $g$ is arbitrary, including
a mixture distribution of certain distribution family, our definition of no community structure includes the case that each weight is generated from a specific distribution in a list of distributions in random order. 
Importantly, when we shuffle the off-diagonal elements $\bb{W}_n$ at random (in element-wise manner), the community structure always disappears. Indeed, in such a case, each element $w_{i, j}'$ of the shuffled matrix $\bb{W}_n'$   independently and identically follows the mixture distribution consisting of different components, i.e.,  $\sum_{k, k'}\pi_{k, k'}g_{k, k'}$  where $\pi_{k, k'}$ is the proportion of elements of cluster block $(k, k')$ for the original matrix $\bb{W}_n$. We use this property for our statistical test as an alternative way for estimating confidential intervals (Section~\ref{stattest}).

\subsection{Statistical test}\label{stattest}
In this section, we develop a statistical test on the existence of community structure defined in Section \ref{setting} 
(i.e., $K=1$ vs. $K>1$). 
We base our test on the asymptotic behavior of the eigenvalues of $S(\bb{W}_n)$ as the number of nodes $n$ goes to $\infty$ when there is no community structure.  A useful result of Random Matrix Theory in our context is that if the elements of an infinite dimensional symmetric matrix $\bb{X}$ independently follow a certain distribution with mean zero and variance one, then the empirical (random) distribution of eigenvalue $\lambda$ of $\bb{X_n}/\sqrt{n}$, where $\bb{X_n}$ is the principal submatrix of $\bb{X}$ for the first $n$ rows and columns, converges almost surely to Wigner semicircular distribution as $n$ goes to $\infty$ (semicircular law)
\cite[p.136]{tao2012topics}:
\begin{eqnarray}
  f_{sc}(\lambda) \equiv \frac{1}{2\pi} \sqrt{4-\lambda^2}.
 \label{semi}
\label{semicircular}
\end{eqnarray}
Note that this law holds for any generative distribution of the elements in matrix $\bb{X}$ (as long as independently drawn), which is referred to as universality property of the law. Also, this law holds even if we replace the diagonal elements by zero's.

In order to apply the semicircular law to our context, we consider a normalization mapping of edge-weight matrix $\bb{W}_n$, transforming each element of the matrix as
\begin{eqnarray*}
\no  T: && w_{i, j} \rightarrow S(w_{i, j})/\sqrt{n},
\end{eqnarray*}
where $S$ is the standardization mapping in Eq.(\ref{S}).
Now, let us assume that the elements in an edge-weight matrix $\bb{W}_n$ are generated as in Eq.(\ref{settingw}). In this setting, 
 the semicircular law suggests that if the eigenvalues of $T(\bb{W}_n)$ for sufficiently large $n$ do not follow Wigner semicircular distribution in Eq.(\ref{semicircular}), then, there should be some $K$-way community structure in the graph ($K>1$) because of our assumption in Eq.(\ref{settingw})\footnote{Without the assumption in Eq.(\ref{settingw}), this property does not hold. For instance, one can make a scale-free graph where the eigenvalues do not follow the semicircular law \cite{Nagao}.}.  
 However, the converse argument does not necessarily hold. That is to say, the fact that the eigenvalues of $T(\bb{W}_n)$ follow Wigner semicircular distribution does not imply that there is no community structure (i.e., $K=1$).  A counter example is given as follows (proof in Appendix~\ref{semicircular1}). 
\begin{exmp} \label{theory counter}
Let $\bb{W}_n$ be a $n \times n$ symmetric edge-weight matrix that has $K$-way community structure with the same cluster size ($n/K$) as defined in Section~\ref{setting}.  Suppose that 
$\mu_{k, k'}=0$ for $\forall k, k'$,  $\sigma_{k, k'}^2=0$ for $k \neq k'$, and $\sigma_{k, k}^2=1$. Then, the empirical eigenvalue distribution of $T(\bb{W}_n)$ almost surely converges to Wigner semicircular distribution in Eq.(\ref{semicircular}) as
$n$ goes to $\infty$.  
\end{exmp}

Nonetheless, in our setting, we can show that an additional condition on the eigenvalue distribution for an exponentially mapped edge-weight matrix ensures that the converse argument also holds. For this purpose, we introduce the exponential mapping $Exp$ that transforms each element of $\bb{W}_n$ as
\begin{eqnarray}
\no  Exp:&& w_{i, j} \rightarrow \exp(t\times w_{i, j})~~\mbox{for}~ i \neq j \\
         && w_{i, i}  \rightarrow 0,
\label{et}
\end{eqnarray}
where $t \in \mathbb{R}$ is a tuning parameter (we do not explicitly denote the dependence of $Exp$ on $t$ because of cluttering). Subsequently, we define the normalization mapping $T_e$ for the exponentially
transformed matrix as 
\begin{eqnarray}
T_e: && w_{i, j} \rightarrow S(Exp(w_{i, j}))/\sqrt{n} .
 \label{et2}
\end{eqnarray}

Now, the following theorem provides a necessary and sufficient condition for the existence of community structure (proof in Appendix~\ref{semicircular2}).
\begin{theorem} \label{theory1}
Let $\bb{W}_n$ be a $n \times n$ weight matrix defined in Section~\ref{setting} with the fixed proportion of 
cluster sizes $(r_1, \ldots, r_K)$ and the pairs of parameters $\{ (\mu_{k, k'}, \sigma_{k, k'})\} (k, k'=1, \ldots, K)$. Suppose that there exists the moment generating function $M(t)$ in an open interval containing zero for $g$ ($g$ is defined 
in Eq.(\ref{settingw})). 
Then, the following statements (C1) and (C2) are equivalent:
\begin{itemize}
\item [] (C1) There is no community structure (i.e., $K=1$)
\item [] (C2) The empirical eigenvalue distribution of the following two matrices almost surely 
converges to Wigner semicircular distribution in Eq.(\ref{semicircular}) as $n$ goes to $\infty$ :
$T(\bb{W}_n)$ and $T_e(\bb{W}_n)$ for some real value $t_0 \neq 0$.  
\end{itemize}
\end{theorem}
Theorem~\ref{theory1} \footnote{Alternatively, one may replace the exponential mapping by 
the square mapping, which requires less assumption on the existence of moments for $g$. 
However, the square transformation seems to have less power when we establish a statistical test (as in the following paragraphs), possibly because it is not one-to-one mapping. This observation motivates us to work on the exponential mapping.} 
motivates us to use the semicircular law to establish a statistical test on 
the null hypothesis $H_0$:
\begin{eqnarray} \label{H0}
  H_0: \mbox{There is no community structure}.
\end{eqnarray}

As has been implied in the proof of Theorem~\ref{theory1}, the violation of the semicircular law for $T(\bb{W}_n)$ is related to differences of means ($\mu_{k, k'}$) among cluster blocks, while the violation of the law for $T_e(\bb{W}_n)$ related to differences of variances ($\sigma_{k, k'}^2$).  Hence, if we take into account the eigenvalues of these two matrices, we can capture differences of the first and second moments of the underlying distributions among cluster blocks. 
Practically, to test the null hypothesis $H_0$,  rather than dealing with the distribution of the whole eigenvalues, 
we focus on extreme values of eigenvalues, because the proof of Theorem~\ref{theory1} suggests that 
the extreme eigenvalues may be closely related to the violation of the semicircular law when there is community structure. The largest eigenvalue may be positively deviated from the expected value $2$, or the smallest eigenvalue may be negatively deviated from the expected value $-2$. 
Note that strictly speaking, the independent assumption on weights is broken if we transform them by $T$ or $T_e$ using
empirical mean and standard deviation $\mu_{emp}$ and $\sigma_{emp}$. For simplicity, however,  we ignore such an effect here. 

The behavior of the largest eigenvalue has been well studied in literature when elements of edge-weight matrix $\bb{W}_n$ are independently generated by certain symmetric distribution $g$ (typically Gaussian, otherwise, its density function may be even with less heavier tails than Gaussian distribution) with mean zero and variance one for non-diagonal elements and with mean zero and variance two for diagonal elements. In this setting,
 the largest eigenvalue $\lambda_{max}$ asymptotically follows Tracy-Widom distribution for Gaussian orthogonal ensembles with parameter $\beta=1$:
\begin{eqnarray}
  \lim_{n \rightarrow \infty}P(\lambda_{max}\leq 2 + x/n^{2/3})=F_1(x),
  \label{F1}
\end{eqnarray} 
where $F_1(x) \equiv \exp \{-(1/2)\int_{x}^{\infty} q(y)dy\} (F_2(x))^{1/2}$ with
$F_2(x) \equiv \exp \{- \int_{x}^{\infty}(y-x)q^2(y)dy\}$ where $q(x)$ is the solution of Painlev\'e I\hspace{-.1em}I equation $d^2q/dx^2=xq + 2q^3$ with the boundary condition $q(x)\sim \mbox{Ai}(x)$ as $x \rightarrow \infty$ \cite{tracy1996orthogonal, tracy2009distributions}. Note that Tracy-Widom distribution is for the maximum eigenvalue 
of specific type of symmetric matrix (e.g., Gaussian ensembles) while the semicircular law holds for the distribution of eigenvalues in general type of symmetric matrix (Wigner ensembles). Moreover,
in our framework, the diagonal elements are all zero, which is in a slightly different situation than the conventional assumption for Tracy-Widom distribution. Nevertheless, 
because of the universality property of Tracy-Widom distribution \cite[Theorem 21.4.3,]{ben2011wigner}, 
we can safely apply Eq.(\ref{F1}) to our context (obviously, our context satisfies the condition of universality that the diagonal part should be symmetric with a sub-Gaussian tail). 

Using Tracy-Widom distribution in Eq.(\ref{F1}), we set confidence intervals for our statistical test as follows.
For the normalized edge-weight matrix $T(\bb{W}_n)$, we set the confidence interval $CI_{max}$ of the largest eigenvalue 
$\lambda_{max}$ at level $\alpha$. Since the violation of the semicircular law occurs as the positive deviation from the expected value, we consider the one-sided confidence interval as $(-\infty, q)$ where $q$ is a critical value at significant level $\alpha$, i.e., $P(\lambda_{max}\geq q | H_0)=\alpha$, which is estimated by $F_1(x)$ in Eq.(\ref{F1}) (refer to the shape of its first derivative in Fig.\ref{F2}b). 
If the generative distribution $g$ may not be symmetric or it may be heavy-tailed, one may evaluate the distribution of the largest eigenvalues by means of permutation test for $T(\bb{W}_n)$, though it may require some computation time. 
In addition to the largest eigenvalue, we also consider to test the smallest eigenvalue $\lambda_{min}$, which may cause the violation of the semicircular law (what matters is indeed the largest magnitude of eigenvalue). In this case, the confidence interval 
$CI_{min}$ is given by $(-q, \infty)$. In the similar manner, we consider to test the largest and the smallest eigenvalue of 
the exponentially normalized weight matrix. We first standardize the data and then apply the mapping $T_e$ where we set $t_0$ to 1/2 as default. 
This results in the transformed matrix
$T_e(S(\bb{W}_n))$ (we denote the confidence intervals as $CI'_{max}$ and $CI'_{min}$, respectively).
Since this procedure involves series of four statistical tests, we set the level of significance to $\alpha/4$ for each test, taking into account the Bonferroni correction (Algorithm~\ref{algo} \footnote{$\mathbb{I}(a)$ is an indicator function: 
1 for correct $a$; 0 otherwise. }).
\begin{algorithm}
\caption{Testing on the existence of community structure}
\label{algo}
\begin{algorithmic} 
\STATE {\bfseries Input:} Edge-weight matrix $\bb{W}$, confidential intervals $CI_{max}$, $CI_{min}$, $CI'_{max}$ and $CI'_{min}$
at level $\alpha/4$.
\STATE $s \leftarrow 0$
\STATE $s$ $\leftarrow$ $s$+$\mathbb{I}$ (max. eigenvalue of $T(\bb{W})$ $\in$ $CI_{max}$)
\STATE $s$ $\leftarrow$ $s$+$\mathbb{I}$ (min. eigenvalues of $T(\bb{W})$ $\in$ $CI_{min}$)
\STATE $s$ $\leftarrow$ $s$+$\mathbb{I}$ (max. eigenvalue of $T_e(S(\bb{W}))$ $\in$ $CI'_{max}$)
\STATE $s$ $\leftarrow$ $s$+$\mathbb{I}$ (min. eigenvalue of $T_e(S(\bb{W}))$ $\in$ $CI'_{min}$)
\IF {$s=4$}
\STATE Accept $H_0$
\ELSE
\STATE Reject $H_0$
\ENDIF
\end{algorithmic}
\end{algorithm}

\section{Simulation study on synthetic data}\label{shynthetic}
In this section, we report on a simulation study to evaluate the performance of our method. First, we investigate the 
validity of using $F_1(x)$ in Eq.(\ref{F1}) for approximation of the distribution of the maximum eigenvalue $\lambda_{max}$ when $n$ is finite. Second, we investigate the power of our method when the null hypothesis $H_0$ is not true. 

Third, we compare the performance of our method outlined in Algorithm~\ref{algo} with other methods. Basically, existing methods in literature consist of two steps. In the first step, a clustering solution for a given graph is yielded by a (arbitrary) clustering method. The yielded solution is subsequently compared with those clustering solutions for randomized graphs, which is further evaluated by a particular statistic. In this study, we adapt one of the state-of-the-art methods based on clustering entropy (\lq CE\rq, originally designed for a unweighted graph) \cite{gfeller2005finding}:
$S = - \frac{1}{L} \sum_{(i, j)} \{ p_{i,j} \log_2p_{i, j} +  (1-p_{i, j}) \log_2(1-p_{i,j}) \}$
where $L$ is the total number of edges in the graph, and $p_{i, j}$ is \lq in-cluster probability\rq~that measures the proportion of accordance of cluster memberships of nodes $i$ and $j$ between the given graph and the randomized graphs
over a number of different noisy contaminations (we set the number of such contaminations to 100). 
As regards clustering, to the best of our knowledge, there is no clustering method that is specifically designed to detect community structure based on differences of patterns of distributions. As a bail-out procedure,  
we consider one of the state-of-the-art methods for signed networks: Signed spectral clustering based on
normalized signed Laplacian method (\lq SignedSpec\rq), which is designed to detect the weakly balanced structure of graph, i.e., positive weights within a cluster and negative weights between clusters \cite{kunegis2010spectral}. We also consider the conventional spectral clustering (normalized Laplacian method, \lq ConvSpec\rq), which is applicable to a graph with positive weights. 
To apply the method \lq ConvSpec\rq~to our context, we transform an edge-weight matrix into a positively-weighted matrix by subtracting   $\underset{i, j}{\mbox{min}}~w_{i, j}$ from each weight.  
Note that the method   \lq ConvSpec\rq~is 
equivalent to the method \lq SignedSpec\rq~when edge weights are all positive. 

\begin{figure}[ht!]
     \begin{center}
        \subfigure[No community]{
            \label{alpha}
            \includegraphics[scale=0.13, trim=15mm 0mm 15mm 0mm]{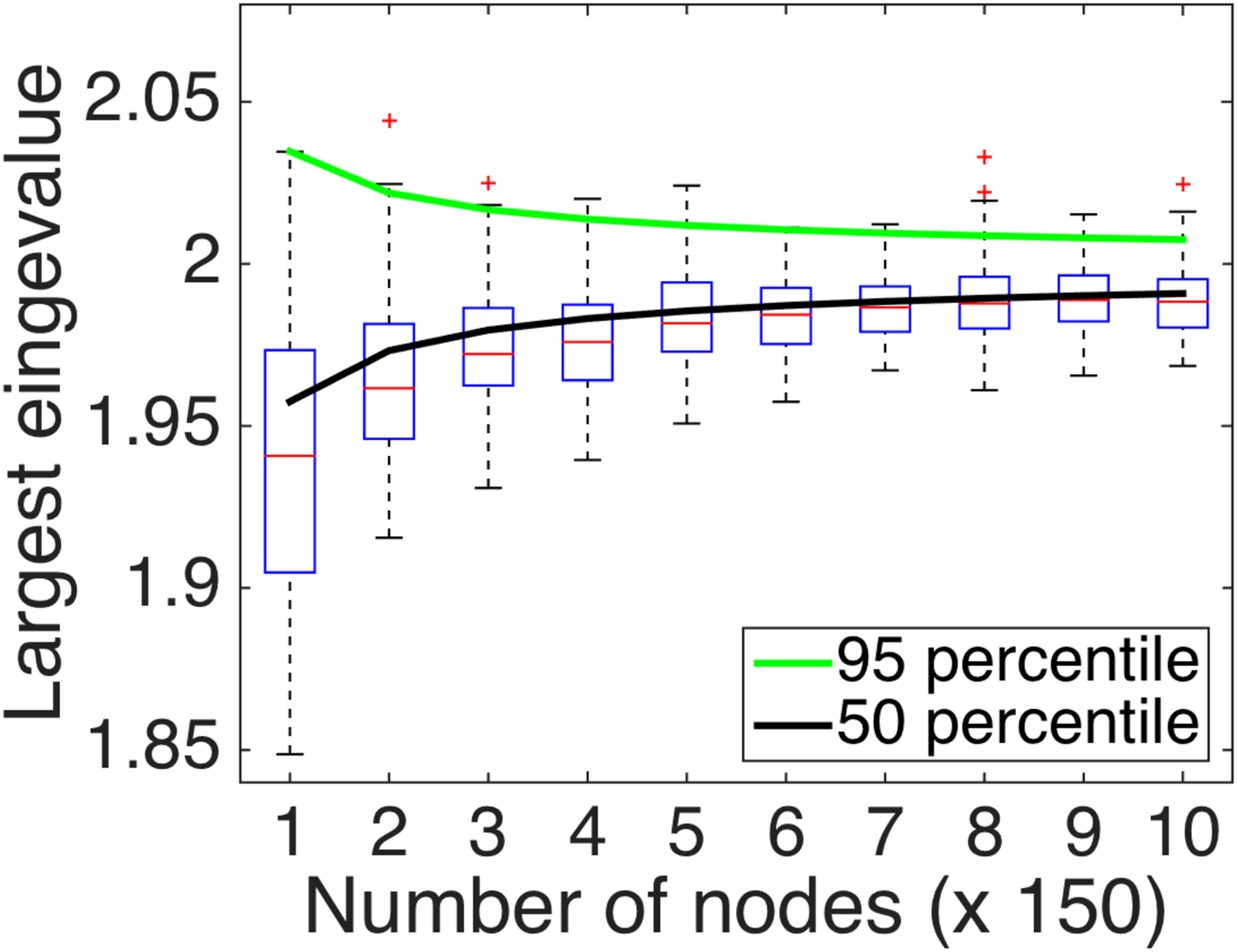} 
        }
        \subfigure[$\mu$ differs]{
           \label{power1}
            \includegraphics[scale=0.13, trim=15mm 0mm 15mm 0mm]{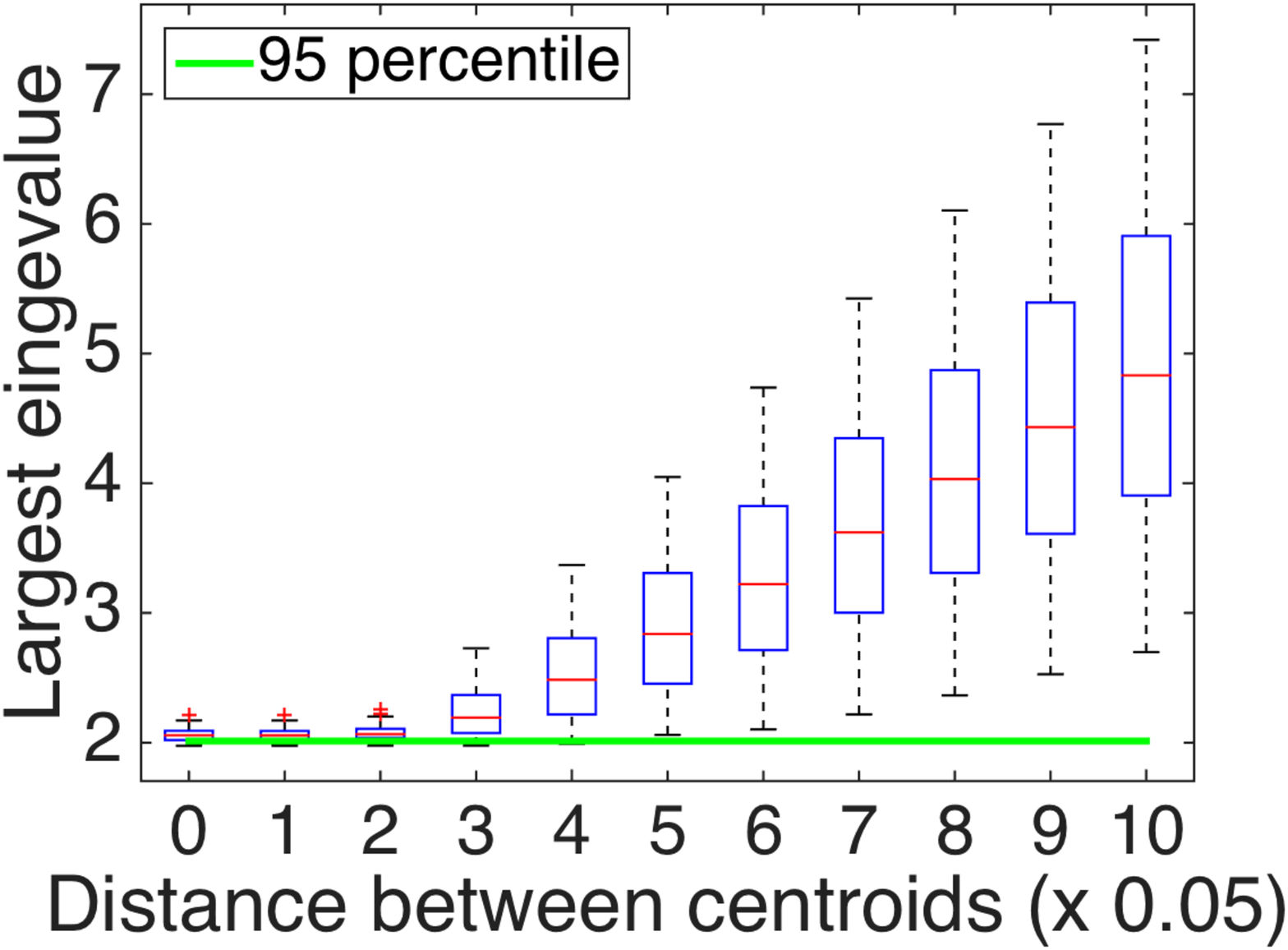} 
        } 
        \subfigure[$\sigma$ differs]{
            \label{power2}
           \includegraphics[scale=0.13, trim=15mm 0mm 15mm 0mm]{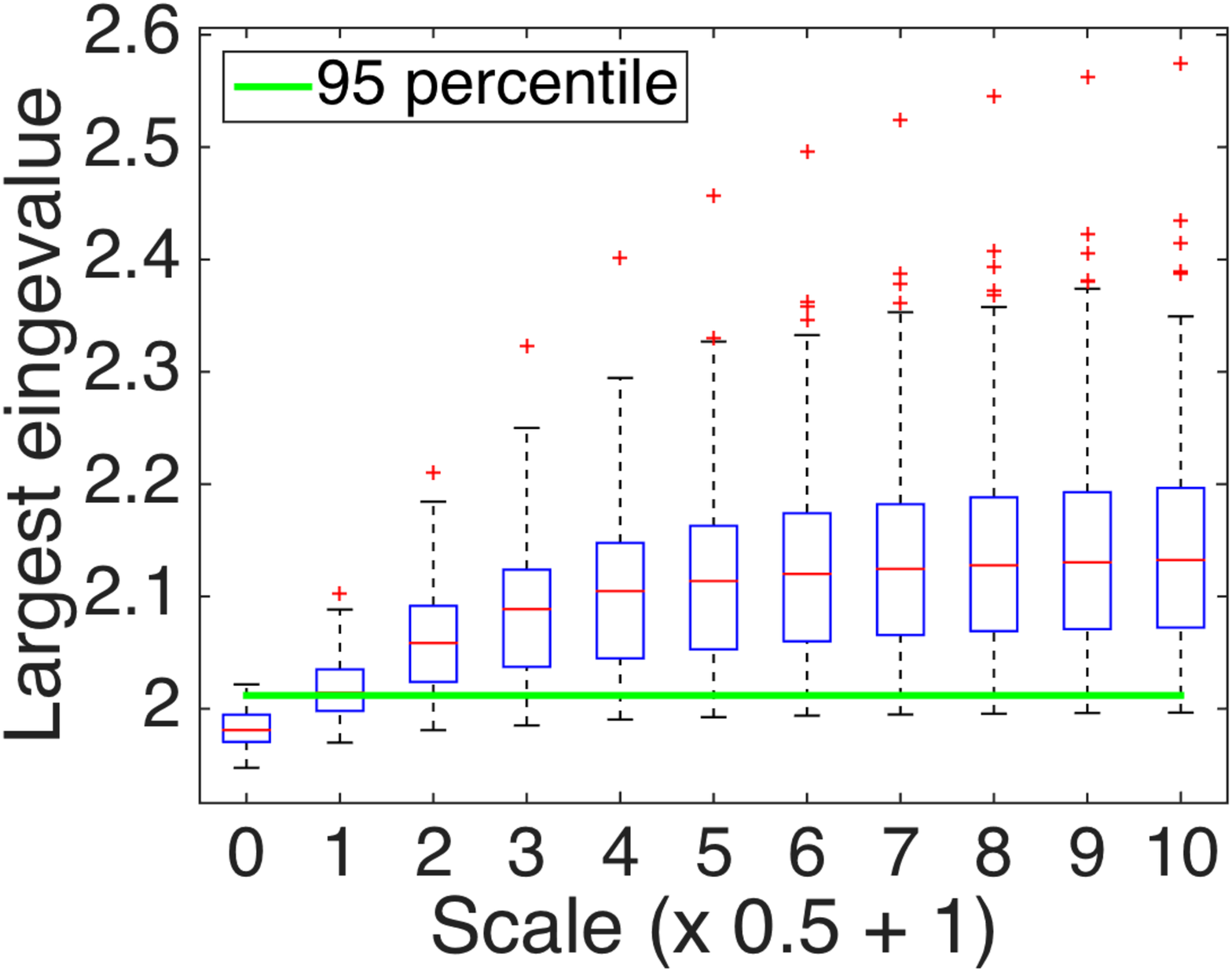} 
        }
        \subfigure[$\sigma$ differs ({\it Exp})]{
           \label{power3}
             \includegraphics[scale=0.13, trim=15mm 0mm 15mm 0mm]{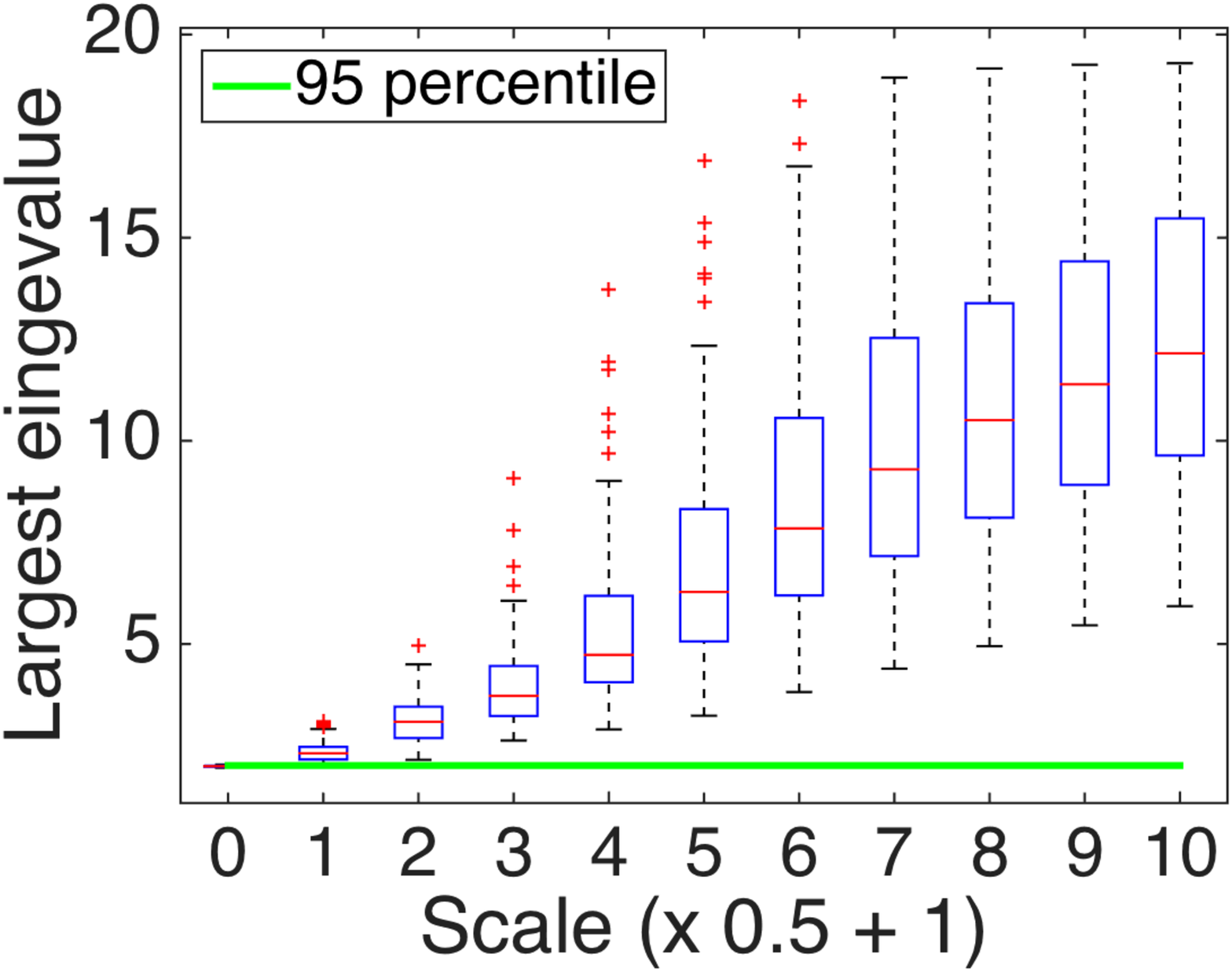}  
        }\\ 
       \end{center}
    \vspace{-5mm}
        \caption{\it \small Boxplots represent 
        distributions of largest eigenvalues for various settings. Panel (a): No community case ($K=1$) of Gaussian ensembles for different number of nodes from 150 to 1500 in x-axis. 
        Panel (b): Five-way community case with the number of nodes 750 and  cluster size (50, 100, 150, 200, 250). Each cluster block is 
        characterized by means of Gaussian distribution (while fixing variance=1), which is randomly chosen from $\{-\mu, \mu\} $ with equal probabilities. The value of $\mu$ is manipulated from 0 to 0.5 of width 0.1 in x-axis.  Panel (c) : Five-way community 
        case characterized by variance (while fixing mean=0), which is randomly chosen from $\{ 1, \sigma^2\} $ with equal probabilities. The value of $\sigma$ is manipulated from 1 to 6 
       of width 1 in x-axis. Panel (d): Five-way community case in the same setting as in (c), but each edge-weight
       matrix is transformed by the exponential mapping {\it Exp} in Eq.(\ref{et}) with $t_0=1/2$. In all panels, the green line denotes the 95 percentile of the largest eigenvalue under the null hypothesis $H_0$ in~(\ref{H0}).}
   \label{check}
\end{figure}
\begin{figure}[ht!]
     \begin{center}
        \subfigure[$\mu$ differs]{
            \label{compar1}
          \includegraphics[scale=0.22, trim=0mm 0mm 0mm 0mm]{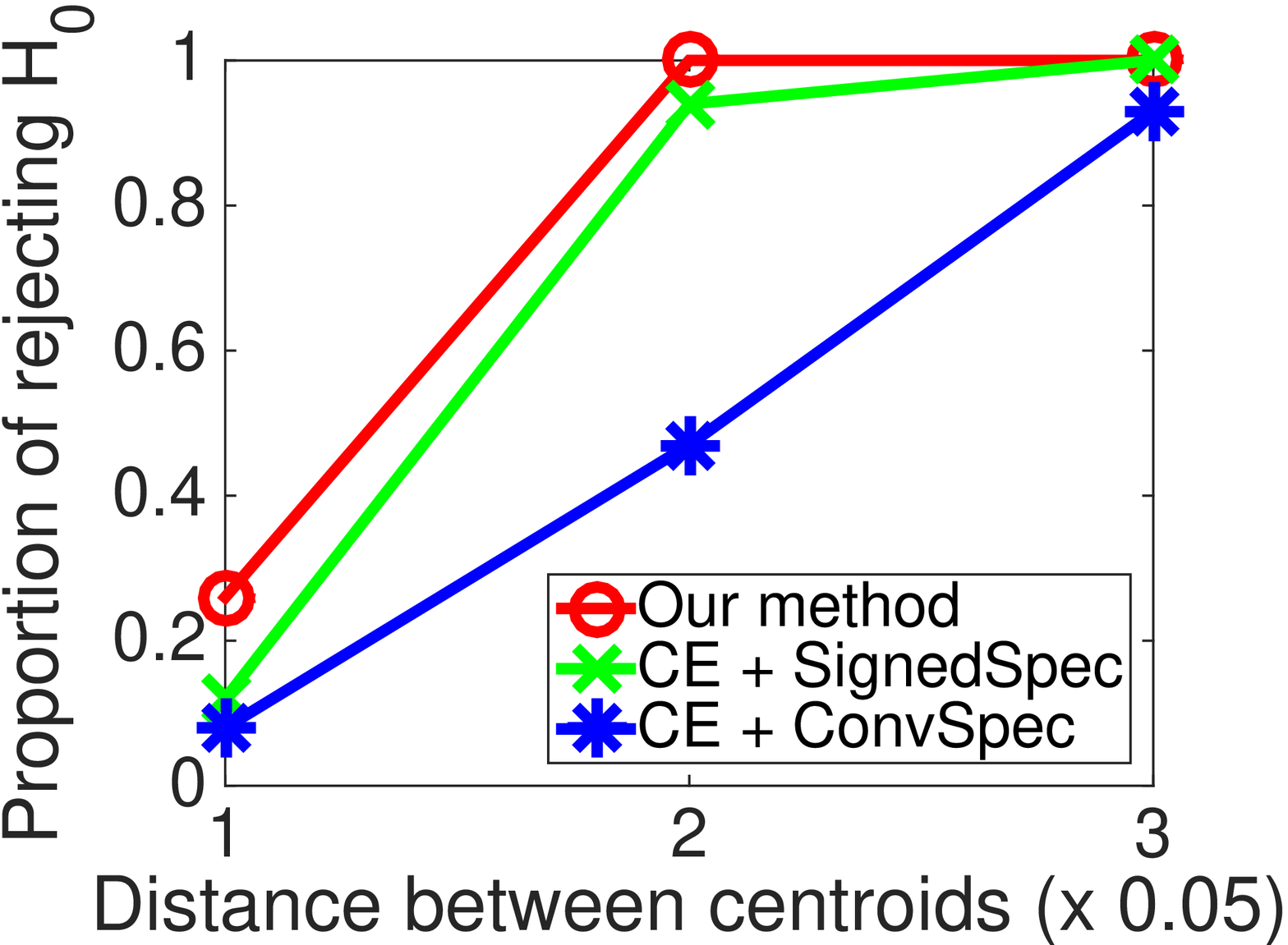} 
           }
        \subfigure[$\sigma$ differs]{
           \label{power}
           \includegraphics[scale=0.22, trim=0mm 0mm 0mm 0mm]{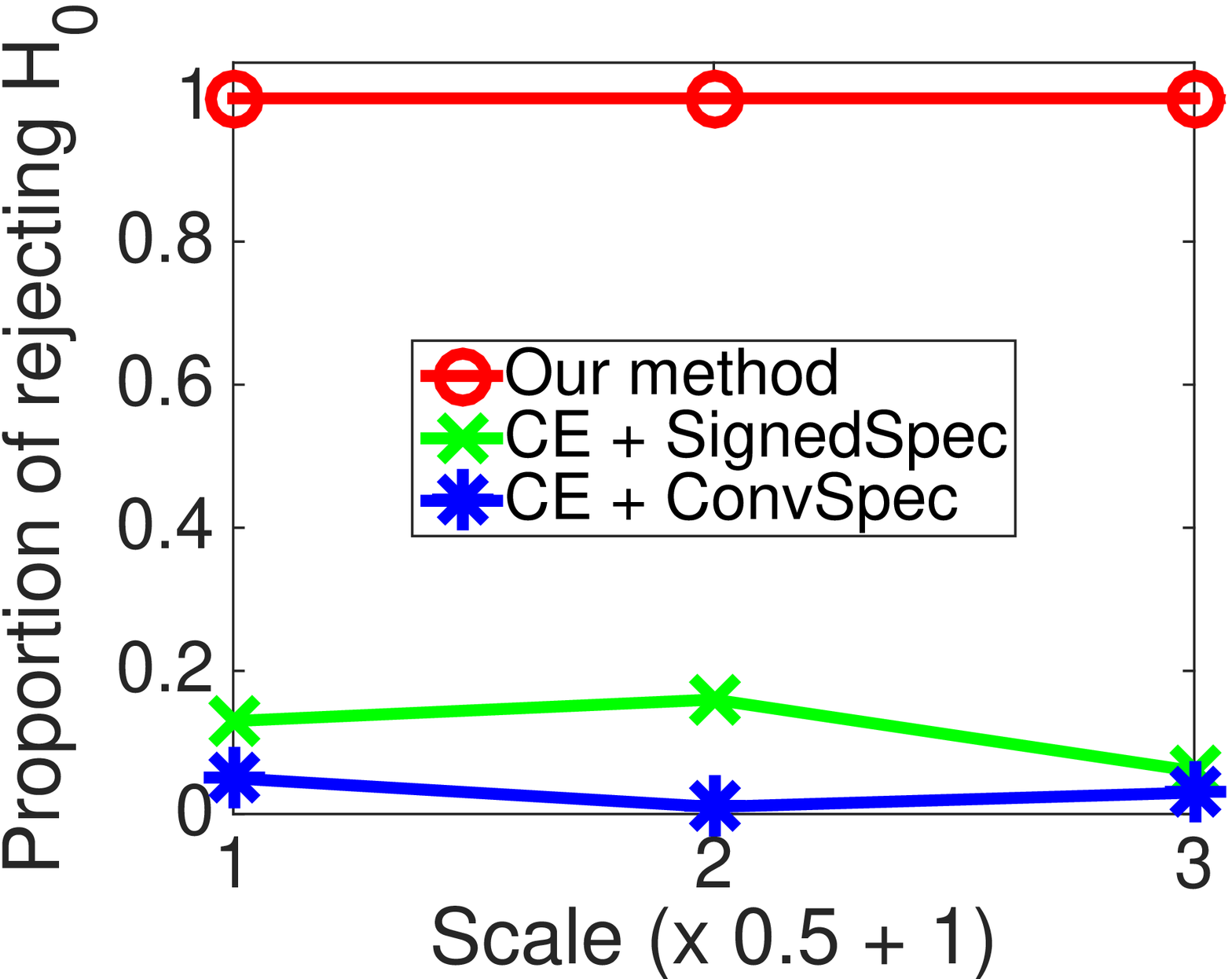} 
        } \\ 
    \end{center}
    \vspace{-5mm}
        \caption{\it \small Comparison of power of test for three different methods: our method, clustering entropy method for yielded cluster solution by signed spectral clustering method (CE + SignedSpec), and clustering entropy method by conventional spectral clustering (CE + ConvSpec). The true community structure is set as follows: cluster size (50, 100, 150, 200, 250);
the means and the variances are manipulated in x-axis of Panel (a) and (b) as in Fig.\ref{check}b and 
Fig.\ref{check}c,  respectively.
 }
   \label{compar2}
\end{figure}

\subsection{Data generation}
For the data structure in this simulation study, we adopted that in \cite{hsieh2012low}, setting the number of clusters to five 
and cluster size to $(10s, 20s, 30s, 40s, 50s)$ where we manipulated integer $s$. In this setting, we have $5 \times 5 = 25$ cluster blocks. In each cluster block, weights were independently drawn from a Gaussian distribution $N(\mu_{k, k'}, \sigma_{k, k'})$ where $\mu_{k, k'}$ and $\sigma_{k, k'}^2$ are the mean and the variance for a cluster block $(k, k')$.
We generated 100 datasets for each setting. 
 
\subsection{Results} \label{simresults}
When the number of nodes ranges from 150 to 1500,  the distribution function $F_1(x)$ in Eq.(\ref{F1}) provides 
a good approximation for the critical value at significance level of $\alpha = 0.05$ under the null hypothesis $H_0$ (Fig.\ref{check}a). Since the function $F_1(x)$ provides the asymptotic probability distribution, this result suggests that 
the function $F_1(x)$ also provides a good approximation for the critical value when the number of nodes goes up more than this range. As regards statistical power, it is implied that our method can well detect the existence of community structure when means $\mu_{k ,k'}$ in each block are different at most by 0.3 
($3 \times 0.05 + 3 \times 0.05 $) when $\sigma_{k, k'}=1$  with the number of nodes being 750 (Fig.\ref{check}b). On the other hand, the power may not be sufficient when differences among cluster blocks are characterized 
by variances $\sigma_{k, k'}^2$ (Fig.\ref{check}c). However, the application of our method to the exponentially transformed matrix by $Exp$ considerably improves the power (Fig.\ref{check}d). All these results suggest good performance of our method to test the existence of  community structure in a graph.

Lastly, we compare the performance of our method with the remainder of the methods. We applied our method as outlined in Algorithm~\ref{algo} to the synthetic data, setting $\alpha$ to 0.05 (hence, $\alpha/4 = 0.0125$). When the community structure is characterized by mean differences, the performance of our method is comparable with the clustering entropy method with signed spectral clustering (CE + SignedSpec), while it outperforms the clustering entropy method with conventional spectral clustering (CE + ConvSpec) (Fig.\ref{compar2}a). On the other hand, when the community structures is characterized by scale differences,  our method considerably outperforms the remainder of the methods (Fig.\ref{compar2}b).

\section{Application to real data} \label{application}
In this section, we experiment our method to real data. The objective is to evaluate the performance of our method when it is applied to various types of real graph data. 

\subsection{Data}
First, we applied our method to the following benchmark graph datasets:
Karate club, {\it Karate} \cite{zachary1977information}; co-authorships in network science, 
{\it Co-authours} \cite{newman2006finding}; Tribal relationships in highland New Guinea, {\it Gahuku-Gama} \cite{read1954cultures}.
The datasets of {\it Karate} and {\it Co-authours} are binary (i.e., $\{0, 1\}$), while the edges in the dataset of {\it Gahuku-Gama} take discrete signed values, $\{-1, 0, 1\}$. The number of nodes for these datasets are 34, 1589, and 16, respectively. These datasets have been well studied in terms of detecting community structure \cite{yang2007community}. The 
clustering results and the underlying social relationships between subjects (nodes) were fully discussed in literature,  clearly suggesting the existence of community structure in these datasets.

\begin{figure}[ht!]
     \begin{center}
           \includegraphics[scale=0.22, trim=0mm 0mm 0mm 0mm]{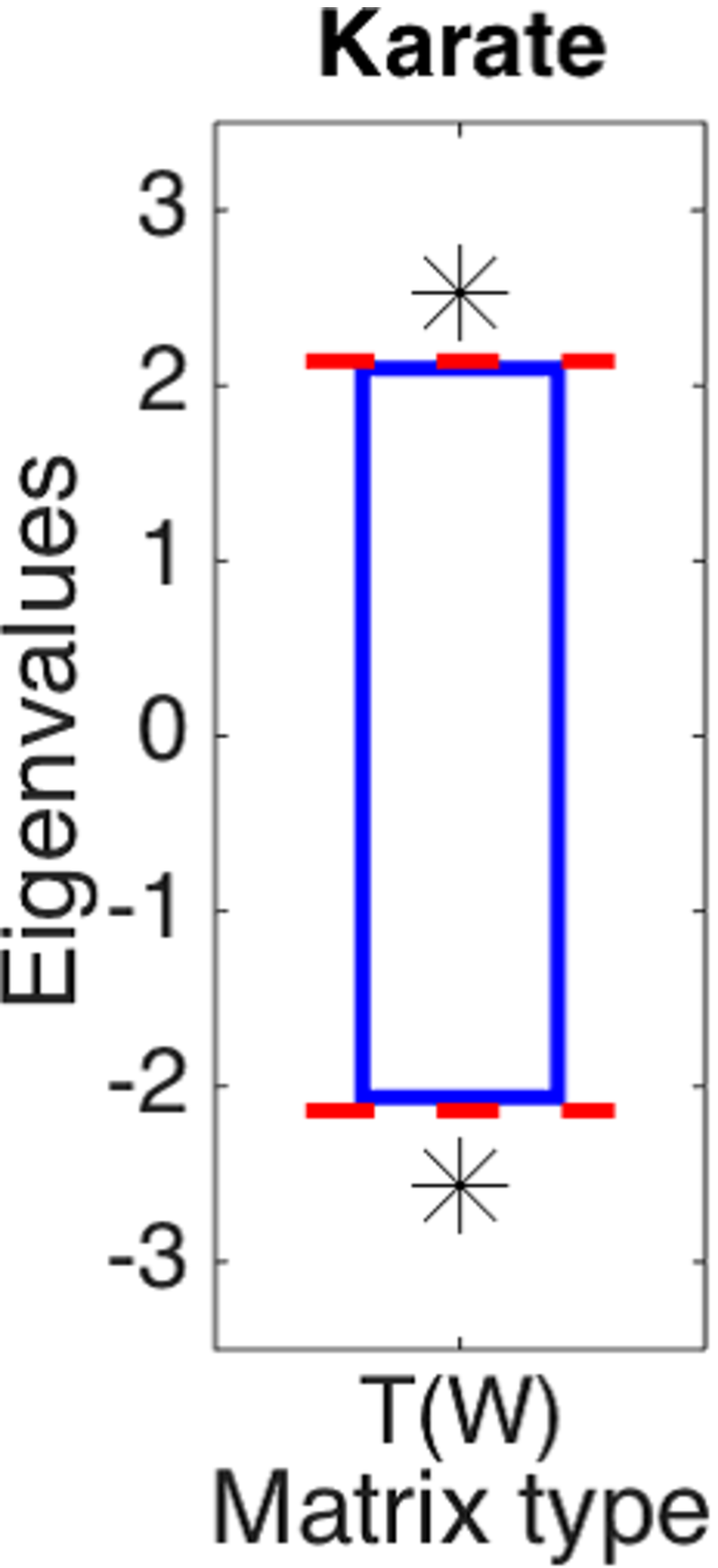} 
           \includegraphics[scale=0.22, trim=0mm 0mm 0mm 0mm]{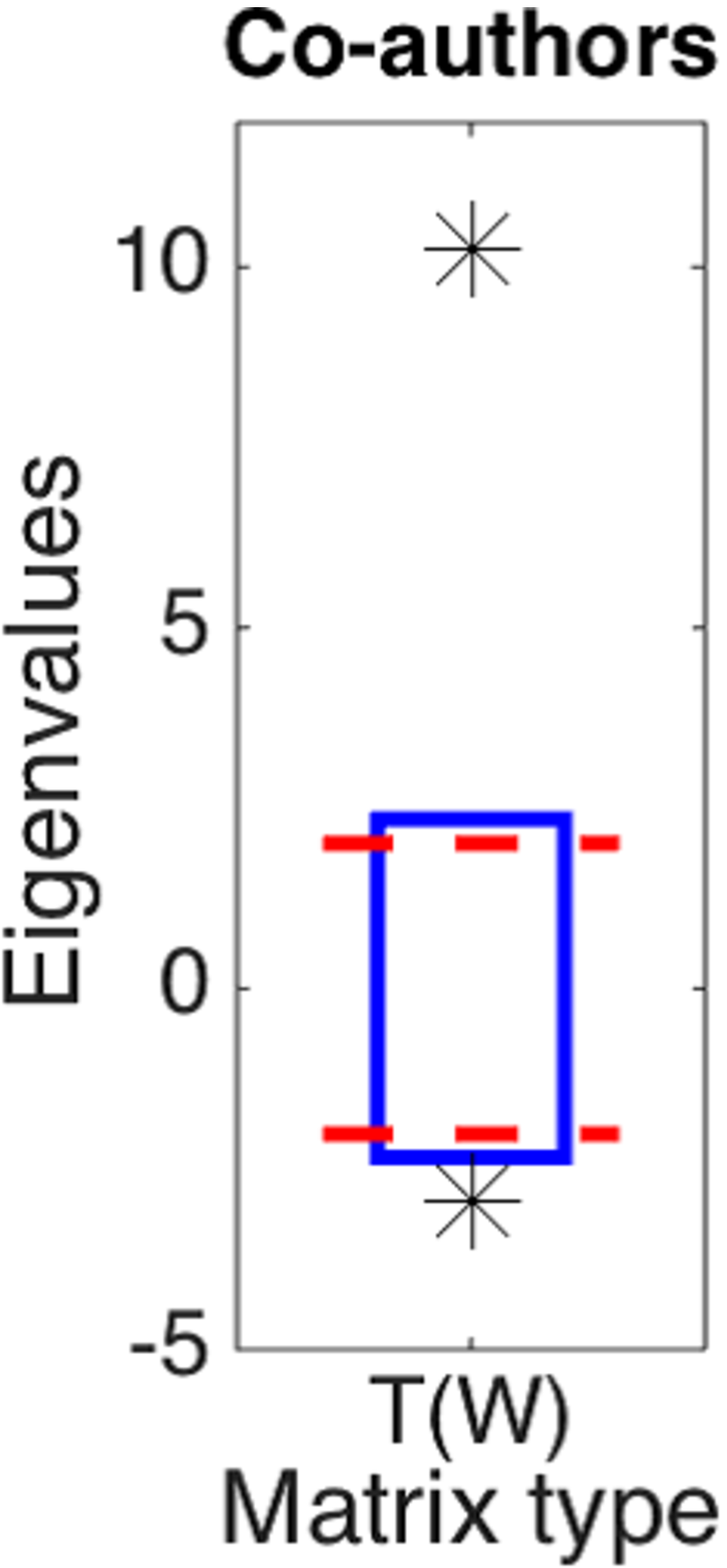}          
          \includegraphics[scale=0.225, trim=0mm 5mm 0mm 0mm]{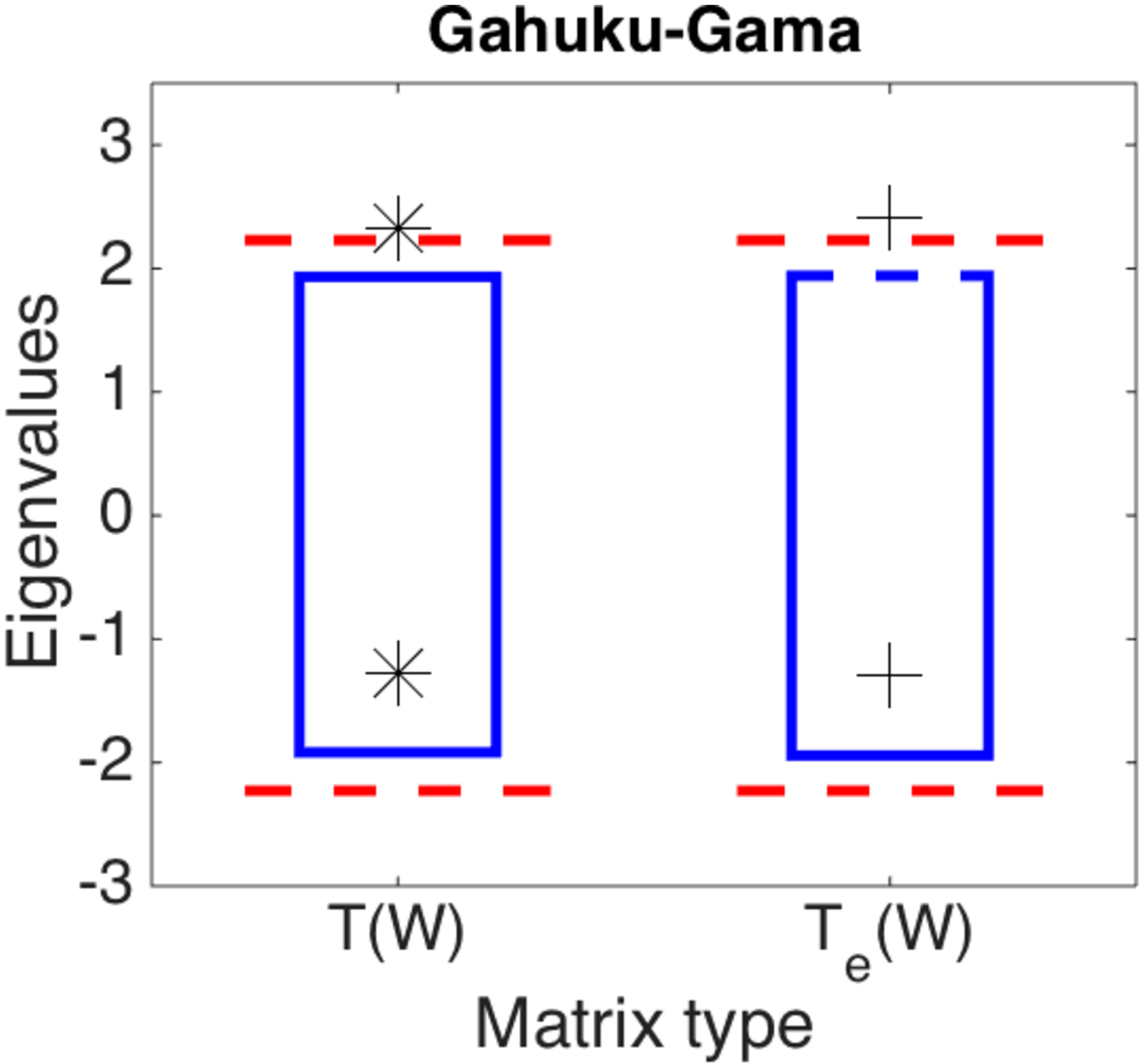}  
     \end{center}
        \caption{\it \small Results of application of our method to real datasets, {\it Karate}, {\it Co-authors}, and {\it Gahuku-Gama} from left to right panels.
   The star maker denotes the maximum or the minimum eigenvalues of the normalized matrix 
        $T(\bb{W})$, while the cross marker denotes those of the exponentially normalized matrix  $T_e(S(\bb{W}))$.
        The edges on top or bottom of boxes denote critical values of these eigenvalues 
        at significance level $\alpha /2$ with $\alpha$=0.05 for {\it Karate} and {\it Co-authors} datasets, and 
        $\alpha/4$ for {\it Gahuku-Gama} dataset. These critical values were yielded by permutation test with 1000 randomized realizations. 
        In contrast, the red dashed lines denote the critical values derived from Tracy-Widom distribution $F_1(x)$. }
   \label{realdata}
\end{figure}
 
Second, we applied our method to a real-valued edge-weighted graph: resting state functional 
MRI ({\it fMRI}) data \cite{fmri}.
The original dataset consists of the level of BOLD (Blood-Oxygen-Level Dependent) signal at short intervals, which reflects neural activity at each tiny portion of the brain, called \lq voxel\rq~(4949 voxels in this dataset). We pre-processed this dataset by evaluating temporal correlations among these voxels and carrying out Fisher's z-transformation for them, which results in a $4949 \time 4949$ edge-weight matrix $\bb{W}$. The objective in this dataset is to experiment our method to such a real-valued edge-weight matrix and to draw a useful implication from the analysis.

\subsection{Results}
\begin{figure}[ht!]
     \begin{center}
           \includegraphics[scale=0.2, trim=0mm 0mm 0mm 0mm]{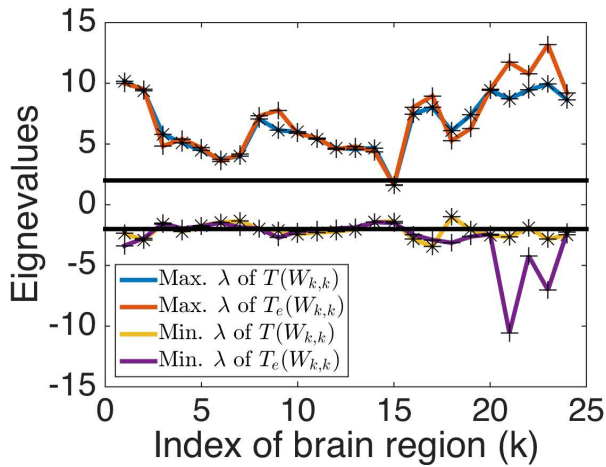} 
        \end{center}
        \caption{\it \small Results of application of our method to the {\it fMRI} dataset.
        The star markers denote the maximum or the minimum of eigenvalues $\lambda$ for normalized weight matrices by the mapping $T$ in various brains regions  with edge-weight matrix $\bb{W}_{k, k}$, indexed by the brain region $k$ in x-axis.
        The cross markers denote those counterparts for exponentially 
        normalized weight matrices by the mapping $T_e$. The horizontal lines denote lines $y=-2$ and $y=2$, which correspond the values to which the minimum and the maximum eigenvalues asymptotically converge. 
        }
   \label{fmriall}
\end{figure}

For the first group of real datasets, our method correctly suggests that the community structure may exist (i.e., $K>1$), whether we estimate critical values either by Tracy-Widom distribution or by permutation test (Fig.\ref{realdata}). 
Note that in the binary case, we always have the same results of our test for the original matrix and the exponentially transposed matrix, because $T(\bb{W})=T_e(S(\bb{W}))$. So, we carried out our test only for $T(\bb{W})$ in {\it Karate} and {\it Co-authors} datasets, setting the significance level to $\alpha/2$. 

For {\it fMRI} dataset, our test rejected the null hypothesis $H_0$, yielding 
the maximum and the minimum eigenvalues as 31.0 and -7.2 for $T(\bb{W})$, and 
31.8 and -10.9 for $T_e(S(\bb{W}))$, which provides strong evidence that there exists the community structure 
in this graph. Furthermore, we carried out our test for subsets of voxels in brain regions that are anatomically predefined, where the number of voxels ranges from 13 to 498. Our test results suggest that there may exist the community structure in each region (except for brain region 16) (Fig.\ref{fmriall}). This result supports the conjecture on heterogeneity of brain activities in anatomically defined brain regions discussed in the literature of neuroscience \cite{birn2001spatial}.

\section{Discussion}
We have proposed a novel method for statistical test on the existence of community structure in an undirected graph that 
is characterized by the first and the second moments of generative model for edge weights.  This method can be considered as an (nontrivial) extension of the recently proposed method \cite{bickel2016hypothesis} from a binary-valued graph to a real-valued one. 
Unlike the existing methods for real-valued graphs, our method does not need a cluster solution. Hence, we can apply this method even to the nontrivial case of clustering in which edge weights take both positive and negative real values. Also, in our approach, we can avoid a nontrivial problem of determining the number of clusters. 
Further, our method is quite efficient in terms of computation time: We need only to evaluate the eigenvalues of edge-weight matrix just once if we use Tracy-Widom distribution, which is due to the asymptotic results provided by Random Matrix Theory. 

As the next step of analysis, one may wonder how to find community memberships when our test rejects the null hypothesis of $K=1$. The present paper did not address this issue, but, it would be quite useful to examine eigenvectors of the edge-weight matrix as in the case of spectral clustering. It is conjectured that some of the eigenvectors of $T(\bb{W})$ and $T_e(S(\bb{W}))$ may have information on community memberships. It would be an important future research topic on how to determine and synthesize relevant eigenvectors for inferring the underlying community structure.
\vspace{2mm}

\newpage
\appendix
\begin{center}
      {\bf \Large APPENDIX}
\end{center}
\section{Proof of Example~\ref{theory counter}} \label{semicircular1}
The edge-weight matrix can be represented as
\begin{eqnarray*}
\bb{W}_n = 
\begin{pmatrix}
  \bb{W}_{n, 1, 1} & \ldots & \bb{0}_{n_1, n_K} \\
  \ldots & \ldots & \ldots \\
  \bb{0}_{n_K,  n_1} & \ldots & \bb{W}_{n, K, K}
\end{pmatrix}
,
\end{eqnarray*}
where $\bb{W}_{n, k, k}$ is the principal submatrix of $\bb{W}_n$ that consists of cluster block $(k, k)$; $n_k$ the number of nodes in the $k$th cluster; $ \bb{0}_{n_k, n_{k'}}$ a $n_k \times n_{k'}$ zero matrix ($1 \leq k, k' \leq K$). Because of the assumption that $\mu_{k, k'}=0$, the normalized matrix $T(\bb{W}_n)$ also has zero off-diagonal blocks. 
As a result, the eigenvalues of $T(\bb{W}_n)$ consist of those eigenvalues of cluster block $(k, k)$ in $T(\bb{W}_n)$. So, it suffices to show that the eigenvalues of these cluster blocks follow the semicircular law. 

Since the variance $\sigma_{k, k}^2$ for $\bb{W}_{n, k, k}$ is one, it becomes that the variance of all elements in $\bb{W}_{n}$ is $\sum_{k=1}^K n_k^2/n^2$. Further, because of the assumption that $n_k  = n/K$ ($k=1, \ldots, K$), 
the variance of all elements in $\bb{W}_n$ becomes $1/K$. So, the standardized matrix $S(\bb{W}_n)$ is given by 
\begin{eqnarray*}
S(\bb{W}_n) = \sqrt{K}\times \bb{W}_n.
\end{eqnarray*}
Hence, the normalized matrix becomes
\begin{eqnarray}
  T(\bb{W}_n) =(\sqrt{K}/\sqrt{n})\times \bb{W}_n.
  \label{coef}
\end{eqnarray}
Since $n =K \times n_k$,  the coefficient $\sqrt{K}/\sqrt{n}$ in Eq.(\ref{coef}) becomes $1/\sqrt{n_k}$. Therefore, the eigenvalues that are relevant for cluster block $(k, k)$ are those eigenvalues of the matrix $\bb{W}_{n, k, k}/\sqrt{n_k}$. 
This suggests that the distribution of these eigenvalues follows the semicircular law as $n$ goes 
to $\infty$. This completes the proof. 

\section{Proof of Theorem~\ref{theory1}} \label{semicircular2}
In our setting, if there is no community structure,  then, the eigenvalues of $T(\bb{W}_n)$ and the eigenvalues of $T_e(\bb{W}_n)$ follow the semicircular law because of universality property of the law. Now, we prove the converse. 
For this purpose, we prove that if there is community structure, then, it violates the semicircular law. 
First, we consider the situation that there is community structure such that some of means $\mu_{k, k'}$ in
 $S(\bb{W}_n)$ differ, while variances $\sigma_{k, k'}$ in  $S(\bb{W}_n)$ are the same across different cluster blocks. 
 Note that means and variances are defined for the standardized matrix $S(\bb{W}_n)$. For simplicity of notation, we denote the normalized weight matrix (i.e., $T(\bb{W}_n)$) as $\bb{W}_n$. 
The matrix $\bb{W}_n$ can be decomposed as follows:
\begin{eqnarray*}
 \bb{W}_n = \bb{W}'_n + \bb{M}_n,
\end{eqnarray*}
where $\bb{M}_n$ is the normalized mean matrix,
\begin{eqnarray*}
  \bb{M}_n = 
  \begin{pmatrix}
      \bb{\mu}_{1, 1} & \ldots & \bb{\mu}_{1, K} \\
      \ldots & \ldots & \ldots \\
      \bb{\mu}_{K, 1} & \ldots & \bb{\mu}_{K, K} 
  \end{pmatrix}
  /\sqrt{n},
\end{eqnarray*}
where $\bb{\mu}_{k, k}=\mu_{k, k} \bb{1}_{n_k}\bb{1}_{n_k}^T$; $n_k$ the number of nodes 
in the $k$th cluster; $\bb{1}_m$ a $m \times 1$ vector with elements one. 
By the dual Weyl inequality \cite[p.40]{tao2012topics} (note that both $\bb{W}'_n$ and $\bb{M}_n$ are symmetric matrices), 
\begin{eqnarray*}
   \lambda_{i+j-n} (\bb{W}_n)\geq \lambda_i(\bb{W}_n') + \lambda_j(\bb{M}_n),
\end{eqnarray*}
where $\lambda_i(\bb{A})$ denotes the $i$th eigenvalue of matrix $\bb{A}$ in descending order.  Letting $i=n$ and $j=1$,
\begin{eqnarray}
   \lambda_{1} (\bb{W}_n) \geq \lambda_n(\bb{W}'_n) + \lambda_1(\bb{M}_n).
   \label{eq2}
\end{eqnarray}
The eigenvalues of $\bb{W}'_n$ follow the semicircular law from its definition. From the semicircular law and Bai-Yin theorem \cite[p.129]{tao2012topics} (because of the existence of the moment-generating function, the fourth moment exists), the eigenvalue $\lambda_n(\bb{W}'_n)$ almost surely converges to  $-2$ as $n \rightarrow \infty$.  On the other hand, we evaluate a lower bound of $\lambda_1(\bb{M}_n)$ as follows. The number of nodes $n_k$ is given by $n_k = r_k \times n$ where $r_k$ is the proportion of the nodes in that cluster. The largest magnitude of eigenvalues  is given 
by the operator norm:
\begin{eqnarray}
\max (|\lambda_1(\bb{M}_n)|, |\lambda_n(\bb{M}_n)|) = \underset{|v|=1}{\mbox{sup}}|\bb{M}_nv|.
\label{operator}
\end{eqnarray}
For simplicity, we assume the left-hand side is given by the largest eigenvalue $\lambda_1(\bb{M}_n)$
(the following argument is applicable when it is $- \lambda_n(\bb{M}_n)$ as well). We evaluate a lower bound of $\lambda_1(\bb{M}_n)$ using Eq.(\ref{operator}). Letting $v=(\bb{v}_1^T, \ldots, \bb{v}_K^T)^T$ with
 $\bb{v}_k = v_k \times \bb{1}_{n_k}$, 
\begin{eqnarray*}
 |\bb{M}_nv|^2 = \frac{1}{n}  \sum_{k=1}^K n_k (\sum_{k'=1}^K n_{k'} \mu_{k, k'}v_{k'})^2
   = n^2 \sum_{k=1}^K r_k  (\sum_{k'=1}^K r_{k'} \mu_{k, k'}v_{k'})^2.
\end{eqnarray*}
Hence, 
\begin{eqnarray}
  \lambda_1(\bb{M}_n) \geq  \epsilon^2 \times n,
  \label{epi}
\end{eqnarray}
where $\epsilon^2 = \sum_{k=1}^K r_k  (\sum_{k'=1}^K r_{k'} \mu_{k, k'}v_{k'})^2$. Because of our assumption, there exists non zero $\mu_{k, k'}$ in Eq.(\ref{epi}). This suggests that by appropriately choosing $\bb{v}$, it becomes that $\epsilon^2 \neq 0$. Therefore, it is concluded that the largest eigenvalue of $\bb{W}_n$ in Eq.(\ref{eq2}) takes (infinitely) larger value than 2 as $n$ goes to $\infty$. This implies the violation of the semicircular that we aim to prove. 

So far, we have assumed  that the variances of cluster block are the same in $\bb{W}$. We relax this condition. Let us assume that variances differ. It suffices to show that $\lambda_n(\bb{W}'_n)$ in Eq.(\ref{eq2}) is lower bounded. We prove  this when $K=2$, but the proof can be easily extended to $K>2$ case.  By the dual Weyl inequality,
\begin{eqnarray*}
 \lambda_n( \bb{W}'_n) \geq \lambda_n(\bb{W}_a') + \lambda_n(\bb{W}_b'),
\end{eqnarray*}
where 
\begin{eqnarray*}
 \bb{W}'_a = 
 \begin{pmatrix}
   \bb{W}'_{11}  & \bb{0} \\
   \bb{0}   & \bb{W}_{22}'
 \end{pmatrix}
 , ~
   \bb{W}'_b =
   \begin{pmatrix}
     \bb{0}  & \bb{W}'_{12} \\
   \bb{W}'_{21}   & \bb{0}     
     \end{pmatrix}
     ,
\end{eqnarray*}
where $\bb{W}'_{i, j}$ is the submatrix of $\bb{W}'$ for cluster block $(i, j)$, and $\bb{W'}^{T}_{12} = \bb{W}'_{21}$.
Since each diagonal block of $\bb{W}'_a$ 
follows the semicircular law up to scale, $\lambda_n(\bb{W}_a')=O(1)$. On the other hand, we augment the diagonal blocks for the symmetric matrix $\bb{W}_b'$ by generating elements from the same distribution as in its off-diagonal blocks. 
\begin{eqnarray*}
  \bb{W}_{b, arg}'=\bb{W}_{b}' + \bb{W}_{b, diag}' ,
 \end{eqnarray*}
where $\bb{W}_{b, arg}'$ is the augmented matrix and  $\bb{W}_{b, diag}'$ the diagonal block matrix. 
Again, by the dual Weyl inequality,
\begin{eqnarray*}
\lambda_{n}(\bb{W}_{b}') \geq \lambda_n(\bb{W}_{b, arg}') +\lambda_n(- \bb{W}_{b, diag}' ).
\end{eqnarray*}
Since the eigenvalues of $\bb{W}_{b, arg}'$ and the eigenvalues of $-\bb{W}_{b, diag}'$ 
follow the semicircular law up to scale, it becomes that
 $\lambda_n(\bb{W}_b')$ is lower bounded. Hence, it is concluded that $ \lambda_n( \bb{W}'_n)$ is also lower bounded. 

Second, we consider the case when means are zero, but some variances differ across cluster blocks. In such a case, 
we consider the exponential transformation of the edge-weight $\bb{W}$ where $w_{i, j} \rightarrow \exp(t_0\times w_{i, j})$.
By definition, the expectation of the new variable $\exp(t_0\times w_{i, j})$ is evaluated at $t=t_0$ by the moment-generating function $M_{k, k'}(t) $ for the distribution of $w_{i, j}$ in block $(k, k')$  where 
$M_{k, k'}(t)$ is given by $\exp(\mu_{k, k'}t)M(\sigma_{k, k'}t)$. In general, the probability distribution is uniquely determined by the corresponding moment-generating function in an open interval containing zero 
\cite[p.155]{rice2006mathematical}. So, if some cluster blocks have different distributions, there exists $t_0$ near zero such that some of $M_{k, k'}(t_0)$ differ for some real $t_0 \neq 0$. This suggests that some means of the new variables differ.  Moreover, if we take $t_0$ small enough such that $M_{k, k'}(4t_0)$ exists (i.e., at least the fourth moment exists for the exponentially transformed variables; indeed, it is possible to do so, because $t_0$ can be taken as much small as one wants), then, the same argument (that we have developed in the first case of different means) can be applied to this case about the violation of the semicircular law (because Bai-Yin theory is also applicable to this case). This completes the proof that if there is community structure, then, the semicircular law is violated.

\bibliographystyle{plain}


\end{document}